\documentclass[a4paper,11pt]{article}
\usepackage{amsmath}
\usepackage{amsfonts}
\usepackage{amssymb}
\usepackage{amscd}
\usepackage{pb-diagram}
\usepackage{epic}

\newtheorem{prop}{Proposition}

\newtheorem{teo}{Theorem}

\date{}
\begin{document}

\title{Corwin-Greenleaf multiplicity function for compact extensions of $\mathbb{R}^n$}
\author{Majdi Ben Halima and Anis Messaoud}
\maketitle
\begin{abstract} Let $G=K\ltimes\mathbb{R}^n$, where $K$ is a compact connected subgroup of $O(n)$
acting on $\mathbb{R}^n$ by rotations. Let
$\mathfrak{g}\supset\mathfrak{k}$ be the respective Lie algebras
of $G$ and $K$, and $pr:
\mathfrak{g}^{*}\longrightarrow\mathfrak{k}^{*}$ the natural
projection. For admissible coadjoint orbits
$\mathcal{O}^{G}\subset\mathfrak{g}^{*}$ and
$\mathcal{O}^{K}\subset\mathfrak{k}^{*}$, we denote by
$n(\mathcal{O}^{G},\mathcal{O}^{K})$ the number of $K$-orbits in
$\mathcal{O}^{G}\cap pr^{-1}(\mathcal{O}^{K})$, which is called
the Corwin-Greenleaf multiplicity function. Let
$\pi\in\widehat{G}$ and $\tau\in\widehat{K}$ be the unitary
representations corresponding, respectively, to $\mathcal{O}^G$
and $\mathcal{O}^K$ by the orbit method. In this paper, we
investigate the relationship between
$n(\mathcal{O}^G,\mathcal{O}^K)$ and the multiplicity
$m(\pi,\tau)$ of $\tau$ in the restriction of $\pi$ to $K$. If
$\pi$ is infinite-dimensional and the associated little group is
connected, we show that $n(\mathcal{O}^G,\mathcal{O}^K)\neq 0$ if
and only if $m(\pi,\tau)\neq 0$. Furthermore, for $K=SO(n)$,
$n\geq 3$, we give a sufficient condition on the representations
$\pi$ and $\tau$ in order that
$n(\mathcal{O}^G,\mathcal{O}^K)=m(\pi,\tau)$.
\end{abstract}
{\footnotesize{\textbf{Mathematics Subject Classification 2010.}
22E20, 22E45, 22E27, 53C30
\vspace{0,2cm}}}
\section{Introduction}
Let $G$ be a connected and simply connected nilpotent Lie group
with Lie algebra $\mathfrak{g}$. It was pointed out by Kirillov
that the unitary dual $\widehat{G}$ of $G$ is parametrized by
$\mathfrak{g}^{*}/G$, the set of coadjoint orbits. The bijection
$$\widehat{G}\simeq\mathfrak{g}^{*}/G$$ is called the Kirillov
correspondence (see [6]). The important feature of this
correspondence is the functoriality relative to inclusion
$K\subset G$ of closed connected subgroups. It means that if we
start with unitary representations $\pi\in\widehat{G}$ and
$\tau\in\widehat{K}$ and if we denote by
$\mathcal{O}^{G}\subset\mathfrak{g}^{*}$ and
$\mathcal{O}^{K}\subset\mathfrak{k}^{*}$ their corresponding
coadjoint orbits, then the multiplicity $m(\pi,\tau)$ of $\tau$ in
the direct integral decomposition of the restriction
$\pi\big\vert_{K}$ can be computed in terms of the space
$\mathcal{O}^{G}\cap p^{-1}(\mathcal{O}^{K})/K$, where $p:
\mathfrak{g}^{*}\longrightarrow\mathfrak{k}^{*}$ denotes the
natural projection. More precisely, Corwin and Greenleaf proved
that the multiplicity $m(\pi,\tau)$ coincides almost everywhere
with the ``mod $K$'' intersection number
$n(\mathcal{O}^{G},\mathcal{O}^{K})$ defined as follows:
$$n(\mathcal{O}^{G},\mathcal{O}^{K}):=\sharp\Big[\big(\mathcal{O}^{G}\cap p^{-1}(\mathcal{O}^{K})\big)/K\Big]$$
(see [2]). The function
$$n:\mathfrak{g}^{*}/G\times\mathfrak{k}^{*}/K\longrightarrow\mathbb{N}\cup\{\infty\},\:\:\:\:\:(\mathcal{O}^{G},\mathcal{O}^{K})\longmapsto
n(\mathcal{O}^{G},\mathcal{O}^{K})$$ is known as the
Corwin-Greenleaf multiplicity function. \vspace{0,2cm}

Suppose now that $G= K\ltimes\mathbb{R}^n$, where $K$ stands for a
compact connected subgroup of the automorphism group
$Aut(\mathbb{R}^n)$.  As usual, $\mathbb{R}^n$ can be equipped
with an Euclidean scalar product which embeds the compact group
$K$ as a subgroup of orthogonal transformations. In the sequel, we
shall simply assume that $K$ is a closed connected subgroup of
$O(n)$. The multiplication law in $G$ is given by
$$(A,a)\cdot(B,b)=(AB,a+Ab)$$
for $(A,a),(B,b)\in G$. As shown by Lipsman in [10], each
irreducible unitary representation of $G$ can be constructed by
holomorphic induction from an admissible linear functional of
$\mathfrak{g}$. Furthermore, two irreducible representations in
$\widehat{G}$ are equivalent if and only if their respective
linear functionals are in the same coadjoint orbit. Thus we can
identify the dual space $\widehat{G}$ with the lattice of
admissible coadjoint orbits.\vspace{0,2cm}

Let $\pi\in\widehat{G}$ and $\tau\in\widehat{K}$ correspond to
admissible coadjoint orbits $\mathcal{O}^{G}$ and
$\mathcal{O}^{K}$ respectively, and let $p:
\mathfrak{g}^{*}\longrightarrow\mathfrak{k}^{*}$ be the
restriction map. In the spirit of the orbit method for nilpotent
Lie groups, one expects that the multiplicity of $\tau$ in
$\pi\big\vert_{K}$ is given by
$\sharp\Big[\big(\mathcal{O}^{G}\cap
p^{-1}(\mathcal{O}^{K})\big)/K\Big]$. Unlike the nilpotent case,
this fact is not in general valid for the semidirect product $G=
K\ltimes\mathbb{R}^n$ (see Theorem 3 below). Assuming that
$\pi\big\vert_{K}$ is multiplicity free as a representation of
$K$, one can ask the following question:
\vspace{0,2cm}\\
\noindent{\bf Question.} Is the intersection $\mathcal{O}^{G}\cap
p^{-1}(\mathcal{O}^{K})$ a single $K$-orbit, provided it is not
empty ?
\vspace{0,2cm}\\
Our interest for this question is motivated by recent
multiplicity-free results in the orbit method obtained by
Kobayashi and Nasrin [9],[14] (see also [11]).\vspace{0,2cm}

By Mackey's little group theory [12,13], the set $\widehat{G}$ is
given by the following procedure. Let $u$ be a non-zero vector in
$\mathbb{R}^n$. We denote by $\chi_{u}$ the unitary character of
$\mathbb{R}^n$ given by $\chi_{u}(v)=e^{iu^{t}v}$ for all
$v\in\mathbb{R}^n$. The stabilizer of $u$ in $K$, denoted by
$K_u$, is called the little group at $u$. For any
$\sigma\in\widehat{K_u}$, define
$\sigma\otimes\chi_u\in\widehat{K_u\ltimes\mathbb{R}^n}$ by
$$(\sigma\otimes\chi_{u})(A,a)=e^{iu^{t}a}\sigma(A)$$
for $A\in K_u$ and $a\in \mathbb{R}^n$. The induced representation
$$\pi_{(\sigma,\chi_{u})}:=Ind_{K_{u}\ltimes
\mathbb{R}^n}^{\,G}(\sigma\otimes\chi_{u})$$ is irreducible and
every infinite-dimensional irreducible unitary representation of
$G$ is equivalent to some $\pi_{(\sigma,\chi_{u})}$. Apart from
these infinite-dimensional unitary representations
$\pi_{(\sigma,\chi_{u})}$, the finite-dimensional unitary
representations of $K$ also yield finite-dimensional unitary
representations of $G$.\vspace{0,2cm}

Let us fix a non-zero vector $u$ in $\mathbb{R}^n$ and assume that
the group $H:=K_u$ is connected. Let $\sigma_{\nu}$ be an
irreducible representation of $H$ with highest weight $\nu$. For
simplicity, write $\pi_{(\nu,u)}$ instead of
$\pi_{(\sigma_{\nu},\chi_{u})}$ and denote by
$\mathcal{O}_{(\nu,u)}^{G}$ the corresponding admissible coadjoint
orbit of $G$. Given an irreducible representation
$\tau_{\lambda}\in\widehat{K}$ with highest weight $\lambda$, we
denote by $\mathcal{O}_{\lambda}^{K}$ the associated admissible
coadjoint orbit of $K$. As above, the multiplicity of
$\tau_{\lambda}$ in the restriction of $\pi_{(\nu,u)}$ to $K$ is
denoted by $m(\pi_{(\nu,u)},\tau_{\lambda})$. The main results of
the present work are\vspace{0,2cm}\\
\textbf{Theorem A.} We have
$$m(\pi_{(\nu,u)},\tau_{\lambda})\neq 0\Leftrightarrow
n(\mathcal{O}^{G}_{(\nu,u)},\mathcal{O}^{K}_{\lambda})\neq 0.$$
\textbf{Theorem B.} Let $(K,H)=(SO(n),SO(n-1))$ with $n\geq 3$.
Assume that $\nu$ and $\lambda$ are strongly dominant weights of
$H$ and $K$, respectively. Then
$$n(\mathcal{O}^{G}_{(\nu,u)},\mathcal{O}^{K}_{\lambda})\leq 1$$
and hence, $m(\pi_{(\nu,u)},\tau_{\lambda})=
n(\mathcal{O}^{G}_{(\nu,u)},\mathcal{O}^{K}_{\lambda}).$

\section{Coadjoint orbits of $K\ltimes\mathbb{R}^n$}
Let $K$ be a closed connected subgroup of $O(n)$ with Lie algebra
$\mathfrak{k}$. The group $K$ acts naturally on $\mathbb{R}^n$ by
rotations, and then one can form the semidirect product
$G=K\ltimes\mathbb{R}^n$. As set $G=K\times\mathbb{R}^n$ and the
group multiplication is given by
$$(A,a)\cdot(A^{'},a^{'})=(AA^{'},a+Aa^{'})$$
for $(A,a),(A^{'},a^{'})\in G$. The Lie algebra of this group is
$\mathfrak{g}=\mathfrak{k}\oplus \mathbb{R}^n$ (as vector space)
with Lie bracket
$$[(U,u),(U^{'},u^{'})]=\big(UU^{'}-U^{'}U,Uu^{'}-U^{'}u\big)$$
for $(U,u),(U^{'},u^{'})\in\mathfrak{g}$. The vector dual space
$\mathfrak{g}^{*}$ of $\mathfrak{g}$ can be identified with
$\mathfrak{k}^{*}\oplus (\mathbb{R}^n)^{*}$. The adjoint action of
$G$ on $\mathfrak{g}$ is expressed by the relation
$$Ad_{{}_{G}}((A,a))(U,u)=(AUA^{t},Au-AUA^{t}a),$$
where $(A,a)\in G$ and $(U,u)\in\mathfrak{g}$. Each linear
functional $F$ on $\mathfrak{g}$ can be identified with an element
$(U,u)\in\mathfrak{g}$ via the natural scalar product
\begin{equation*}
\langle(U,u),(V,v)\rangle:=\frac{1}{2}tr(UV^t)+u^tv,
\end{equation*}
where $(V,v)\in\mathfrak{g}$. It follows that for $(A,a)\in G$,
$(U,u)\in\mathfrak{g}^{*}$ and $(V,v)\in \mathfrak{g}$
\begin{eqnarray*}
\langle Ad_{{}_{G}}^*((A,a))(U,u),(V,v)\rangle&=&\langle
(U,u),Ad_{{}_{G}}((A,a)^{-1})(V,v)\rangle\\
&=&\frac{1}{2}tr((AUA^t)V^t)+(Au)^t(Va)+(Au)^tv.
\end{eqnarray*}
Given two vectors $a$ and $b$ in $\mathbb{R}^n$, there exists a
unique matrix $W_{a,b}$ in $\mathfrak{k}$ such that
\begin{eqnarray*}
\frac{1}{2}tr(W_{a,b}V^t)=b^tVa
\end{eqnarray*}
for all $V\in\mathfrak{k}$. Observe that
$$W_{Aa,Ab}=AW_{a,b}A^{t}$$
for all $A\in K$. Now, we can write
\begin{eqnarray*}
\langle
Ad_{{}_{G}}^*((A,a))(U,u),(V,v)\rangle=\langle(AUA^t+W_{a,Au},Au),(V,v)\rangle,
\end{eqnarray*}
i.e.,
\begin{eqnarray*}
Ad_{{}_{G}}^*((A,a))(U,u)&=&(AUA^t+W_{a,Au},Au).
\end{eqnarray*}
Therefore, the coadjoint orbit $\mathcal{O}^{G}_{(U,u)}$ of $G$
through $(U,u)$ is given by
\begin{eqnarray*}
   \mathcal{O}^{G}_{(U,u)} &=&\{(AUA^t+W_{a,Au},Au);
A\in K, a\in \mathbb{R}^n\}\\
&=&\{(A(U+W)A^t,Au); A\in K, W\in\mathcal{W}_u\},
\end{eqnarray*}
where $\mathcal{W}_u:=\{W_{a,u}; a\in \mathbb{R}^n\} $.
\section{Irreducible unitary representations of $K\ltimes\mathbb{R}^n$}
We keep the notation of Section 2. Let $u$ be a non-zero vector in
$\mathbb{R}^n$. We denote by $\chi_{u}$ the unitary character of
the vector Lie group $\mathbb{R}^n$ given by
$\chi_{u}(v)=e^{iu^{t}v}$ for all $v\in\mathbb{R}^n$. We define
the little group $K_{u}$ at $u$ to be the stabilizer of $u$ in
$K$. Let $\sigma$ be an irreducible unitary representation of
$K_{u}$ on some (finite-dimensional!) Hilbert space $\mathcal{H}$.
The map
$$\sigma\otimes\chi_{u}: (A,a)\longmapsto e^{iu^{t}a}\sigma(A)$$
is a representation of the semidirect product $K_{u}\ltimes
\mathbb{R}^n$. Let $L^{2}(K,\mathcal{H})$ be the completion of the
vector space of all continuous maps $\eta: K\longrightarrow
\mathcal{H}$ with respect to the norm
$$\Vert\eta\Vert=\Big(\int_{K}\Vert\eta(A)\Vert^{2}dA\Big)^{\frac{1}{2}},$$
where $dA$ is a normalized Haar measure on $K$. Define
$L^{2}(K,\mathcal{H})^{\sigma}$ to be the subspace of
$L^{2}(K,\mathcal{H})$ consisting of the maps $\xi$ which satisfy
the covariance condition
$$\xi(AB)=\sigma(B^{t})\xi(A)$$
for $B\in K_{u}$ and $A\in K$. The induced representation
$$\pi_{(\sigma,\chi_{u})}:=Ind_{K_{u}\ltimes
\mathbb{R}^n}^{\,G}(\sigma\otimes\chi_{u})$$ is realized on
$L^{2}(K,\mathcal{H})^{\sigma}$ by
$$\pi_{(\sigma,\chi_{u})}((A,a))\xi(B)=e^{i(Bu)^{t}a}\xi(A^{t}B),$$
where $(A,a)\in G$, $\xi\in L^{2}(K,\mathcal{H})^{\sigma}$ and
$B\in K$. Mackey's theory tells us that the representation
$\pi_{(\sigma,\chi_{u})}$ is irreducible and that every
infinite-dimensional irreducible unitary representation of $G$ is
equivalent to some $\pi_{(\sigma,\chi_{u})}$. Furthermore, two
representations $\pi_{(\sigma,\chi_{u})}$ and
$\pi_{(\sigma^{'},\chi_{u^{'}})}$ are equivalent if and only if
$u$ and $u^{'}$ lie in the same $K$-orbit and the representations
$\sigma$ and $\sigma^{'}$ are equivalent under the identification
of the conjugate subgroups $K_{u}$ and $K_{u^{'}}$. In this way,
we obtain all irreducible representations of $G$ which are not
trivial on the normal subgroup $\mathbb{R}^n$. On the other hand,
every irreducible unitary representation $\tau$ of $K$ extends
trivially to an irreducible representation, also denoted by
$\tau$, of $G$ by $\tau(A,a):=\tau(A)$ for $A\in K$ and $a\in
\mathbb{R}^n$.\vspace{0,2cm}

For $\Omega\in \mathbb{R}^n/K$, let $u$ be any element of $\Omega$
and define $\widehat{G}(\Omega)$ to be the set of all induced
representations $\pi_{(\sigma,\chi_{u})}$ with
$\sigma\in\widehat{K_{u}}$. Then $\widehat{G}$ is the disjoint
union of the $\widehat{G}(\Omega)$, $\Omega$ in $\mathbb{R}^n/K$.
Up to identification, we can write
$$\widehat{G}=\widehat{K}\bigcup
\Big(\bigcup_{\Omega\in\Lambda}\widehat{G}(\Omega)\Big)$$ where
$\Lambda$ is the set of non-trivial $K$-orbits in $\mathbb{R}^n$.
Finally, notice that
$\bigcup_{\Omega\in\Lambda}\widehat{G}(\Omega)$ has full
Plancherel measure in $\widehat{G}$ (see [7]). \vspace{0,2cm}

\section{Admissible coadjoint orbits of $K\ltimes\mathbb{R}^n$}
Let us fix a non-zero vector $u$ in $\mathbb{R}^n$ and assume that
the little group $K_u$ is connected. Take $T_{K_u}$ and $T_{K}$ to
be maximal tori respectively in $K_{u}$ and $K$ such that
$T_{K_u}\subset T_{K}$. Consider an irreducible unitary
representation $\sigma_{\nu}: K_{u}\longrightarrow
U(\mathcal{H}_{\nu})$ with highest weight $\nu$. Then
$$\pi_{(\sigma_{\nu},\chi_{u})}=Ind_{K_{u}\ltimes \mathbb{R}^n}^{\,G}(\sigma_{\nu}\otimes\chi_{u})$$
is an irreducible unitary representation of $G$. To simplify
notation, we shall write $\pi_{(\nu,u)}$ instead of
$\pi_{(\sigma_{\nu},\chi_{u})}$. We fix an irreducible unitary
representation $\tau_{\nu}: K\longrightarrow
U(\mathcal{H}^{'}_{\nu})$ with highest weight $\nu$ and we realize
the representation space $\mathcal{H}_{\nu}$ of $\sigma_{\nu}$ as
the smallest $K_{u}$-invariant subspace of $\mathcal{H}^{'}_{\nu}$
that contains the $\nu$-weight space of
$\mathcal{H}^{'}_{\nu}$.\vspace{0,2cm}

Choose a normalized highest weight vector $w_{\nu}$ in
$\mathcal{H}^{'}_{\nu}$ and define a vector
$U_{\nu}\in\mathfrak{k}$ by the relation
$$\displaystyle\frac{1}{2}tr(U_{\nu}U^t)=-i\langle d\tau_{\nu}(U)w_{\nu},w_{\nu}\rangle$$
for all $U\in\mathfrak{k}$. If we set $\ell_{\nu,u}:=(U_{\nu},u)$,
then we can see that the stabilizer $ G(\ell_{\nu,u})$ of
$\ell_{\nu,u}$ in $G$ is equal to
$G(\ell_{\nu,u})=K(\ell_{\nu,u})\ltimes
\mathbb{R}^n(\ell_{\nu,u})$. Hence, $\ell_{\nu,u}$ is aligned in
the sense of Lipsman (see [10]). A linear functional $\ell\in
\mathfrak{g}^* $ is called admissible, if there exists a unitary
character $\chi$ of the connected component of
 $G(\ell) $, such that $d\chi=i\ell\vert_{\mathfrak{g}(\ell)}$.
Notice that the linear functional $\ell_{\nu,u}$ is admissible and
so, according to Lipsman [10], the representation of $G$ obtained
by holomorphic induction from $\ell_{\nu,u}$ is equivalent to the
representation $ \pi_{(\nu,u)} $. Now, for an irreducible unitary
representation $\tau_{\lambda}$ of $K$ with highest weight
$\lambda$, we take the linear functional
$\ell_{\lambda}:=(U_\lambda,0)$ of $ \mathfrak{g}^* $ which is
clearly aligned and admissible. Hence, the representation of $G$
obtained by holomorphic induction from the linear functional
$\ell_{\lambda}$ is equivalent to the representation
$\tau_\lambda$.\vspace{0,2cm}

We denote by $ \mathcal{O}^{G}_\lambda$ the coadjoint orbit of
$\ell_\lambda$ and by $ \mathcal{O}^{G}_{(\nu,u)} $ the coadjoint
orbit of $\ell_{\nu,u}$. Let $ \mathfrak{g}^\ddagger \subset
\mathfrak{g}^*$ be the union of all the $\mathcal{O}^{G}_{(\nu,u)}
$ and of all the $\mathcal{O}^{G}_{\lambda} $ and denote by
$\mathfrak{g}^\ddagger/G$ the corresponding set in the orbit
space. It  follows now from [10], that $ \mathfrak{g}^\ddagger $
is just the set of all admissible linear functionals of $
\mathfrak{g}$. The result of Lipsman stated in the introduction
gives us a bijection between the space $\mathfrak{g}^\ddagger/G$
of admissible coadjoint orbits and the unitary dual $\widehat{G}$.
\section{Main results}
We continue to use the notation of the previous sections. Fix a
non-zero vector $u$ in $\mathbb{R}^n$ and assume that the Lie
subgroup $H:=K_{u}$ is connected. The corresponding Lie algebra of
$H$ is denoted by $\mathfrak{h}$, i.e.,
$\mathfrak{h}=\{U\in\mathfrak{k};\:Uu=0\}.$ Let $
\pi_{(\nu,u)}\in\widehat{G} $ and $\tau_{\lambda}\in\widehat{K}$
be as before. To these irreducible unitary representations, we
attach respectively the coadjoint orbits
$\mathcal{O}^{G}_{(\nu,u)} $ and $\mathcal{O}^{K}_{\lambda}$. Here
$\mathcal{O}^{K}_{\lambda}$ is the orbit in $\mathfrak{k}^{*}$
through $U_{\lambda}$, i.e.,
$\mathcal{O}^{K}_{\lambda}=Ad_{K}^{*}(K)U_{\lambda}$. If
$p:\mathfrak{g}^{*}\longrightarrow\mathfrak{k}^{*}$ denotes the
canonical projection, then we can define the ``mod $K$'' number
$$n(\mathcal{O}^{G}_{(\nu,u)},\mathcal{O}^{K}_{\lambda}):=\sharp\Big[\big(\mathcal{O}^{G}_{(\nu,u)}\cap p^{-1}(\mathcal{O}^{K}_{\lambda})\big)/K\Big].$$
\begin{prop} Let $H_{\nu}$ be the stabilizer of $U_{\nu}$ in $H$.
Then
$$n(\mathcal{O}^{G}_{(\nu,u)},\mathcal{O}^{K}_{\lambda})=\sharp\Big[\big(\big(U_{\nu}+\mathcal{W}_{u}\big)\cap \mathcal{O}^{K}_{\lambda}\big)/H_{\nu}\Big].$$
\end{prop}
\textbf{Proof.} Assume that $\mathcal{O}^{G}_{(\nu,u)}\cap
p^{-1}(\mathcal{O}_{\lambda}^{K})\neq\emptyset$, i.e., the set
$$\mathcal{F}:=\big\{W\in\mathcal{W}_u;\:U_{\nu}+W\in\mathcal{O}^{K}_{\lambda}\big\}$$
is non-empty. We define an equivalence relation in $\mathcal{F}$
by
$$W_1\sim W_2\Leftrightarrow\exists A\in
H_{\nu};\:W_2=AW_1A^{t}.$$ The set of equivalence classes is
denoted by $\mathcal{F}/H_{\nu}$. Letting
$$\mathcal{E}_W:=\{(A(U_{\nu}+W)A^t,Au); A\in K\}$$
for $W\in\mathcal{F}$, one can easily prove that
$$W_1\sim
W_2\Leftrightarrow\mathcal{E}_{W_1}=\mathcal{E}_{W_2}.$$ Since
$$\mathcal{O}^{G}_{(\nu,u)}\cap
p^{-1}(\mathcal{O}^{K}_{\lambda})=\underset{W\in\mathcal{F}}{\bigcup}\mathcal{E}_{W},$$
we can deduce that
\begin{eqnarray*}
n(\mathcal{O}^{G}_{(\nu,u)},\mathcal{O}^{K}_{\lambda})&=&\sharp\Big[\big(\mathcal{O}^{G}_{(\nu,u)}\cap
p^{-1}(\mathcal{O}^{K}_{\lambda})\big)/K\Big]\\&=&\sharp\Big[\mathcal{F}/H_{\nu}\Big]\\&=&\sharp\Big[\big(\big(U_{\nu}+\mathcal{W}_{u}\big)\cap
\mathcal{O}^{K}_{\lambda}\big)/H_{\nu}\Big].
\end{eqnarray*}
This completes the proof of the proposition.$\hfill\square$\\

Let us denote by $q$ the canonical projection from
$\mathfrak{k}^{*}$ to $\mathfrak{h}^{*}$.
\begin{prop} The intersection $\mathcal{O}^{G}_{(\nu,u)}\cap
p^{-1}(\mathcal{O}^{K}_{\lambda})$ is non-empty if and only if
$U_{\nu}\in q(\mathcal{O}^{K}_{\lambda})$.
\end{prop}
\textbf{Proof.} The result of the proposition immediately follows
from the equivalence $$\mathcal{O}^{G}_{(\nu,u)}\cap
p^{-1}(\mathcal{O}_{\lambda}^{K})\neq\emptyset\Leftrightarrow
\big(U_{\nu}+\mathcal{W}_{u}\big)\cap\mathcal{O}^{K}_{\lambda}\neq\emptyset$$
and the direct sum decomposition
$\mathfrak{k}=\mathfrak{h}\oplus\mathcal{W}_u$.$\hfill\square$

\vspace{0,2cm}

Now we turn our attention to the multiplicity
$m(\pi_{(\nu,u)},\tau_{\lambda})$ of $\tau_{\lambda}$ in the
restriction of $\pi_{(\nu,u)}$ to $K$. We have\begin{eqnarray*}
m(\pi_{(\nu,u)},\tau_{\lambda})&=:&mult(\pi_{(\nu,u)}\big\vert_{K},\tau_{\lambda})\\&=&mult(Ind_{H}^{K}\sigma_{\nu},\tau_{\lambda})\\
&=&mult(\tau_{\lambda}\big\vert_{H},\sigma_{\nu}).
\end{eqnarray*}

\begin{prop} The representation $\tau_{\lambda}$ occurs in the
restriction of $\pi_{(\nu,u)}$ to $K$ if and only if $U_{\nu}\in
q(\mathcal{O}^{K}_{\lambda})$.
\end{prop}
\textbf{Proof.} A proof of the if part can be found in [1]. The
only if part is directly obtained by applying a result of
Guillemin-Sternberg [3,4] (compare [5]) which relates the branching
problem for compact connected Lie groups to the projection of
coadjoint orbits.$\hfill\square$\vspace{0,2cm}

From the above propositions, we immediately obtain
\begin{teo} We have $$m(\pi_{(\nu,u)},\tau_{\lambda})\neq 0\Leftrightarrow n(\mathcal{O}^{G}_{(\nu,u)},\mathcal{O}^{K}_{\lambda})\neq
0.$$
\end{teo}\vspace{0,3cm}

As illustrated by the following example, the multiplicities
$m(\pi_{(\nu,u)},\tau_{\lambda})$ and
$n(\mathcal{O}^{G}_{(\nu,u)},\mathcal{O}^{K}_{\lambda})$ may
coincide in certain cases. Let $1_K$ stand for the trivial
representation of $K$. It is well known that the multiplicity of
$1_{K}\in\widehat{K}$ in $\pi_{(\nu,u)}\big\vert_{K}$ is $0$ or
$1$. The coadjoint orbit corresponding to the representation
$1_{K}$ is $\{0\}$ and the formula for the multiplicity
$n(\mathcal{O}^{G}_{(\nu,u)},\{0\})$ is
$$n(\mathcal{O}^{G}_{(\nu,u)},\{0\})=\sharp\Big[\big(\mathcal{O}^{G}_{(\nu,u)}\cap \mathfrak{k}^{\perp}\big)/K\Big],$$
where
$\mathfrak{k}^{\perp}:=p^{-1}(\{0\})=\{\ell\in\mathfrak{g}^{*};\:\ell(\mathfrak{k})=0\}$.
Clearly, $\mathcal{O}^{G}_{(\nu,u)}\cap \mathfrak{k}^{\perp}$ is a
single $K$-orbit whenever the intersection is non-empty, i.e.,
$n(\mathcal{O}^{G}_{(\nu,u)},\{0\})\leq 1$. Thus,
$m(\pi_{(\nu,u)},\tau_{\lambda})=
n(\mathcal{O}^{G}_{(\nu,u)},\mathcal{O}^{K}_{\lambda}).$
\vspace{0,4cm}

In the remainder of this paper, we fix $K=SO(n)$ with $n\geq 3$.
Then $G=K\ltimes\mathbb{R}^n$ is the so-called Euclidean motion
group. Without loss of generality, we can take
$u=(0,...,0,r)^{t}\in\mathbb{R}^{n}$ with
$r\in\mathbb{R}_{+}^{*}$. In fact, if $u$ and $u^{'}$ belong to
the same sphere centered at zero and of radius $r=\|u\|$, then
$K_{u^{'}}=AK_uA^{t}$ for some $A\in K$ and the representations
$\pi_{(\nu,u)}$ and $\pi_{(\nu,u^{'})}$ are equivalent. Notice
that the little group at the vector
$u=(0,...,0,r)^{t}\in\mathbb{R}^n$ is the subgroup $H=SO(n-1)$. In
this case, $H$ is a multiplicity free subgroup of $K$ in the
following sense: For any irreducible representation $\rho$ of $K$,
all the irreducible $H$-components of the restriction
$\rho\big\vert_{H}$ of $\rho$ to $H$ have multiplicity at most
$1$. The multiplicity free property for the restriction
$\pi_{(\nu,u)}\big\vert_{K}$ may predict that the Corwin-Greenleaf
function $n(\mathcal{O}^{G}_{(\nu,u)},\mathcal{O}^{K}_{\lambda})$
is either $0$ or $1$, and then $m(\pi_{(\nu,u)},\tau_{\lambda})=
n(\mathcal{O}^{G}_{(\nu,u)},\mathcal{O}^{K}_{\lambda})$. Next, we
shall prove that this prediction turns out to be true when the
weights $\nu$ and $\lambda$ are strongly dominant.\vspace{0,2cm}

Let us first recall a useful fact concerning the weight lattice of
$SO(n)$, $n\geq 3$. Fix the Cartan subalgebra of
$\mathfrak{so}(n)$ consisting of the two-by-two diagonal blocks
$$\left(
\begin{array}{cc}
0&\theta_j\\
-\theta_j&0\\
\end{array}\right), j=1,...,\big[\frac{n}{2}\big],$$ starting from the upper left. Here, $[\frac{n}{2}]$ denotes the largest
integer smaller than $\frac{n}{2}$. For an integer
$j\in\{1,...,[\frac{n}{2}]\}$, denote by $e_j$ the associated
evaluation functional on the complexification of the Cartan
subalgebra. If the standard choice of positive roots is made (see,
e.g., [8]), then the dominant weights (resp. strongly dominant
weights) $\lambda$ for $SO(n)$ are given by expressions
\begin{equation*}
\lambda=\lambda_1e_1+...+\lambda_de_d \longleftrightarrow
\lambda=(\lambda_1,...,\lambda_d)
\end{equation*}
such that
$$\lambda_1\geq...\geq\lambda_{d-1}\geq\vert \lambda_d\vert
\:\:\:(\text{resp}.\:\: \lambda_1>...>\lambda_{d-1}>\vert
\lambda_d\vert)$$ when $n=2d$ is even, and
$$\lambda_1\geq...\geq\lambda_d\geq 0\:\:\:(\text{resp}.\:\:
\lambda_1>...>\lambda_d> 0)$$ when $n=2d+1$ is odd, where
$2\lambda_i$ and $\lambda_i-\lambda_j$ are integers for all $i$,
$j$.\vspace{0,2cm}

Now, we are in position to prove
\begin{teo} Let $(K,H)=(SO(n),SO(n-1))$ with $n\geq 3$. Assume that $\nu$ and $\lambda$ are strongly dominant weights of $H$ and $K$, respectively. Then
$$n(\mathcal{O}^{G}_{(\nu,u)},\mathcal{O}^{K}_{\lambda})\leq 1$$
and hence, $m(\pi_{(\nu,u)},\tau_{\lambda})=
n(\mathcal{O}^{G}_{(\nu,u)},\mathcal{O}^{K}_{\lambda}).$
\end{teo}
\textbf{Proof.} We will prove the theorem only for the pair
$(K,H)=(SO(2d+1),SO(2d))$. Analogous proof holds in the remaining
case $(K,H)=(SO(2d+2),SO(2d+1))$.\vspace{0,2cm}

Given the vector $u=(0,...,0,r)^t\in\mathbb{R}^n$,
$r\in\mathbb{R}_{+}^{*}$, we have
\begin{equation*}
\mathcal{W}_{u}= \left\{\left(\begin{array}{cccc}
0 & \ldots & 0 & -y_1\\
\vdots & \ddots & \vdots & \vdots\\
0 & \ldots & 0 & -y_{2d}\\
y_1 & \ldots & y_{2d} & 0\\
\end{array}
 \right);\:y_j\in\mathbb{R}\:\forall\: j\right\}.
\end{equation*}
Letting $\nu=(\nu_1,...,\nu_d)$ be a strongly dominant weight of
$H$ and
$$J= \left(
\begin{array}{cc}
0&1\\
-1&0\\
\end{array}\right),$$ we associate to the
representation $ \pi_{(\nu,u)} $ the linear functional $
\ell_{\nu,u}=(U_{\nu},u)$ in $ \mathfrak{g}^ * $, where
$$U_\nu=
\left( \begin{array}{cccc}
\nu_1J & \ldots & 0&0\\
\vdots & \ddots & \vdots&\vdots\\
0 & \ldots & \nu_dJ&0\\
0 & \ldots & 0&0\\
\end{array}\right).$$
We shall denote by $H_\nu$ the stabilizer of $U_\nu$ in
$H$.\vspace{0,2cm}

Let $\lambda=(\lambda_{1},...,\lambda_{d})$ be a strongly dominant
weight of $K$. We link the representation $\tau_{\lambda}$ to the
linear functional $\ell_{\lambda}=(U_\lambda,0) $ in $
\mathfrak{g}^*$, where
$$U_\lambda=
\left( \begin{array}{cccc}
\lambda_1J & \ldots & 0&0\\
\vdots & \ddots & \vdots& \vdots\\
0 & \ldots & \lambda_dJ&0\\
0 & \ldots & 0&0\\
\end{array}\right).$$ \vspace{0,2cm}
Assume that
$n(\mathcal{O}^{G}_{(\nu,u)},\mathcal{O}^{K}_{\lambda})\neq 0$.
Then there exists a skew-symmetric matrix
\begin{equation*} W=
\left(\begin{array}{cccc}
0 & \ldots & 0 & -y_1\\
\vdots & \ddots & \vdots & \vdots\\
0 & \ldots & 0 & -y_{2d}\\
y_1 & \ldots & y_{2d} & 0\\
\end{array}
 \right)
 \end{equation*}
in $\mathcal{W}_{u}$ such that $U_{\nu}+W=AU_{\lambda}A^{t}$ for
some $A\in K$. For all $x\in \mathbb{R}$, we have
 $$det(U_{\nu}+W-ix\mathbb{I})=i(-1)^{d+1}xP(x)$$ where
$P$ is the unitary polynomial of degree $2d$ given by
\begin{eqnarray*}
P(x)=\prod_{i=1}^d(x^2-\nu_i^2)-\sum_{j=1}^d\Big((y_{2j-1}^2+y_{2j}^2)\prod_{i=1
,i\not=j}^d(x^2-\nu_i^2)\Big).
\end{eqnarray*}
Let $Q$ be the following unitary polynomial of degree $2d$:
$$Q(x)=\prod_{i=1}^d(x^2-\nu_i^2)-\sum_{j=1}^d\dfrac{\prod_{i=1}^{i=d}(\lambda_i^2-\nu_j^2)\prod_{i=1,i\not=j}^{i=d}(x^2-\nu_i^2)}{\prod_{i=1,
i\not=j}^{i=d}(\nu_i^2-\nu_j^2)}.$$ Applying the Lagrange's
interpolation theorem, we have
$$\prod_{i=1}^d(\lambda_{k}^2-\nu_i^2)=
\sum_{j=1}^d\dfrac{\prod_{i=1}^{i=d}(\lambda_i^2-\nu_j^2)\prod_{i=1,i\not=j}^{i=d}(\lambda_{k}^2-\nu_i^2)}{\prod_{i=1,
i\not=j}^{i=d}(\nu_i^2-\nu_j^2)},$$ and so $Q(\pm\lambda_k)=0$ for
$k=1,...,d$. It follows that $P=Q$, and then
\begin{eqnarray*}
y_{2j-1}^2+y_{2j}^2=\dfrac{\prod_{i=1}^{i=d}(\lambda_i^2-\nu_j^2)}{\prod_{i=1,
i\not=j}^{i=d}(\nu_i^2-\nu_j^2)}
\end{eqnarray*}
for all $j=1,...,d$.\vspace{0,2cm}

Consider again the set
$\mathcal{F}=\big\{W\in\mathcal{W}_u;\:U_{\nu}+W\in\mathcal{O}^{K}_{\lambda}\big\}$.
Since $\mathcal{F}$ is stable under the natural action of $H_\nu$
on $\mathcal{W}_u$, the stabilizer $H_\nu$ is necessarily included
in the torus
\begin{equation*}
T_H= \left\{\left(
\begin{array}{cccccccc}
cos(\theta_1)       & sin(\theta_1) &          &\\
-sin(\theta_1) & cos(\theta_1)     &          &\\
        &        &\ddots    &\\
        &        &          &  cos(\theta_d)      & sin(\theta_d) &          &\\
         &       &          &  -sin(\theta_d) & cos(\theta_d)    &          &\\
         &        &         &          &        &          &1\\
\end{array}
\right);\:\theta_j\in\mathbb{R}\:\forall\: j\right\}.
\end{equation*}
Furthermore, it is clear that $T_H\subseteq H_\nu$ and then we
have $H_\nu=T_H$. By observing that the $H_\nu$-action on
$\mathcal{F}$ is transitive, we deduce that
$\mathcal{O}^{G}_{(\nu,u)}\cap p^{-1}(\mathcal{O}_{\lambda}^{K})$
is a single $K$-orbit, i.e.,
$$n(\mathcal{O}^{G}_{(\nu,u)},\mathcal{O}^{K}_{\lambda})=m(\pi_{(\nu,u)},\tau_{\lambda})=1.$$
$\hfill\square$\\

Concluding this section, let us prove the following result:
\begin{teo} Let the pair $(K,H)$ be either $(SO(2d+1),SO(2d))$ or $(SO(2d+2),SO(2d+1))$. If the dominant
weight $\nu=(\nu_1,...,\nu_d)$ of $H$ satisfies
$\nu_1=...=\nu_d=\alpha$ for some
$\alpha\in\frac{1}{2}\mathbb{N}^{*}$, then for any dominant weight
$\lambda$ of $K$ with $\lambda\neq\nu$ we have
$$n(\mathcal{O}^{G}_{(\nu,u)},\mathcal{O}^{K}_{\lambda})\neq 1.$$
Consequently, if
$n(\mathcal{O}^{G}_{(\nu,u)},\mathcal{O}^{K}_{\lambda})\neq 0$
then $m(\pi_{(\nu,u)},\tau_{\lambda})\neq
n(\mathcal{O}^{G}_{(\nu,u)},\mathcal{O}^{K}_{\lambda}).$
\end{teo}
\textbf{Proof.} We take $(K,H)=(SO(2d+1),SO(2d))$. The proof of
the remaining case $(K,H)=(SO(2d+2),SO(2d+1))$ goes along the same
lines.\vspace{0,2cm}

Let $\nu=(\nu_1,...,\nu_d)$ be a dominant weight of $H$ such that
$\nu_1=...=\nu_d=\alpha$ with
$\alpha\in\frac{1}{2}\mathbb{N}^{*}$. Assume that
$n(\mathcal{O}^{G}_{(\nu,u)},\mathcal{O}^{K}_{\lambda})\neq 0$ for
some dominant weight $\lambda$ of $K$. Then there exists a
skew-symmetric matrix
\begin{equation*} W=
\left(\begin{array}{cccc}
0 & \ldots & 0 & -y_1\\
\vdots & \ddots & \vdots & \vdots\\
0 & \ldots & 0 & -y_{2d}\\
y_1 & \ldots & y_{2d} & 0\\
\end{array}
 \right)
 \end{equation*}
in $\mathcal{W}_{u}$ such that $U_{\nu}+W=AU_{\lambda}A^{t}$ for
some $A\in K$. For all $x\in \mathbb{R}$, we have
 $$det(U_{\nu}+W-ix\mathbb{I})=i(-1)^{d+1}xP(x)$$ with
$$
P(x)=(x^2-\alpha^2)^{d-1}\Big(x^2-\alpha^2-\sum_{j=1}^d
(y_{2j-1}^2+y_{2j}^2)\Big).$$ Applying the branching rule from
$SO(2d+1)$ to $SO(2d)$ (see, e.g., [8]), we observe that the
weight $\lambda$ is of the form
$$\lambda=(\beta,\underset{d-1}{\underbrace{\alpha,...,\alpha}})$$ where $\beta\in\frac{1}{2}\mathbb{N}^{*}$ and
$\beta-\alpha\in\mathbb{N}^{*}$. Thus we get
$$\sum_{j=1}^d(y_{2j-1}^2+y_{2j}^2)=r^2$$
with $r=\sqrt{\beta^2-\alpha^2}.$ Since $H_\nu$ is a proper
subgroup of $SO(2d)$, the $H_\nu$-action on the
$(2d-1)$-dimensional sphere centered at zero and with radius $r$
is not transitive. That is, the $H_\nu$-action on the set
$\mathcal{F}=\big\{W\in\mathcal{W}_u;\:U_{\nu}+W\in\mathcal{O}^{K}_{\lambda}\big\}$
is not transitive. Therefore,
$n(\mathcal{O}^{G}_{(\nu,u)},\mathcal{O}^{K}_{\lambda})\neq 1$ and
so
$$m(\pi_{(\nu,u)},\tau_{\lambda})\neq n(\mathcal{O}^{G}_{(\nu,u)},\mathcal{O}^{K}_{\lambda}).$$
$\hfill\square$
\small{

\vspace{0,5cm}

Department of Mathematics, Faculty of Sciences at Sfax, University of Sfax, Route de Soukra, B. P. 1171, 3000-Sfax, Tunisia\\
E-mail adress: majdi.benhalima@yahoo.fr\\

Department of Mathematics, Preparatory Institute for Engineering Studies at Gafsa, University of Gafsa, El Khayzorane street-Zaroug, 2112-Gafsa, Tunisia\\
E-mail adresses: anis.messaoud@ipeig.rnu.tn


\begin{thebibliography}{99}
\vspace{0,2cm}
\bibitem{Arnal} D. Arnal, M. Ben Ammar, M. Selmi, Le probl\`eme de la r\'eduction \`a un sous-groupe dans la quantification par d\'eformation, \emph{Ann. Fac. Sci. Toulouse}, {\bf 12} (1991), 7-27.
\vspace{0,2cm}
\bibitem{Corwin} L. Corwin, F. Greenleaf, Spectrum and multiplicities for unitary representations in nilpotent Lie groups, \emph{Pacific J. Math.}, {\bf 135} (1988), 233-267.
\vspace{0,2cm}
\bibitem{Guillemin} V. Guillemin, S. Sternberg, Convexity properties of the moment mapping, \emph{Invent. math.}, {\bf 67} (1982), 491-513.
\vspace{0,2cm}
\bibitem{Sternberg} V. Guillemin, S. Sternberg, Geometric quantization and multiplicities of group representations,
\emph{Invent. math.}, {\bf 67} (1982), 515-538.
\vspace{0,2cm}
\bibitem{Heckman} G.J. Heckman, Projection of orbits and asymptotic behavior of multiplicities for compact connected Lie groups, \emph{Invent. math.}, {\bf 67} (1982), 333-356.
\vspace{0,2cm}
\bibitem{Kirillov} A. A. Kirillov, \emph{Lectures on the orbit method}, Amer. Math. Soc., Providence, RI, 2004.
\vspace{0,2cm}
\bibitem{Kleppner} A. Kleppner, R.L. Lipsman, The Plancherel formula for group extensions, \emph{Ann. Sci. Ecole Norm. Sup.}, {\bf 4} (1972), 459-516.
\vspace{0,2cm}
\bibitem{Knapp} A.W. Knapp, \emph{Lie Groups Beyond an Introduction\,-\,Second Edition}, Birkh\"auser, Boston, 2002.
\vspace{0,2cm}
\bibitem{Kobayashi} T. Kobayashi and S. Nasrin, Multiplicity one theorem in the orbit method, in: \emph{Lie groups and symmetric spaces}, 161-169, Amer. Math. Soc. Transl. Ser. 2, 210, Amer. Math. Soc., Providence, RI, 2003.
\vspace{0,2cm}
\bibitem{Lipsman} R.L. Lipsman, Orbit theory and harmonic analysis on Lie groups with co-compact nilradical, \emph{J. Math. pures et appl.}, {\bf 59} (1980), 337-374.
\vspace{0,2cm}
\bibitem{Ronald} R.L. Lipsman, Attributes and applications of the Corwin-Greenleaf multiplicity function, \emph{Contemp. Math.}, {\bf 177} (1994), 27-46.
\vspace{0,2cm}
\bibitem{Mackey1} G.W. Mackey, \emph{The theory of unitary group representations}, Chicago University Press, 1976.
\vspace{0,2cm}
\bibitem{Mackey2} G.W. Mackey, \emph{Unitary group representations in physics, probability and number theory}, Benjamin-Cummings, 1978.
\vspace{0,2cm}
\bibitem{Nasrin} S. Nasrin, Corwin-Greenleaf multiplicity functions for Hermitian symmetric spaces and multiplicity-one theorem in the orbit method, \emph{Int. J. Math.}, {\bf 21} (2010), 279-296.

\end{thebibliography}
\end{document}